\documentclass[12pt]{article}

\usepackage[a4paper,margin=2.5cm]{geometry}
\usepackage{amsmath,amsfonts,hyperref}

\title{Sparsity and non-Euclidean embeddings}

\author{Omer Friedland \and Olivier Gu\'edon}

\newcommand\address{\noindent\leavevmode\noindent
Omer Friedland, \\
Institut de Math\'ematiques de Jussieu, \\
Universit\'e Pierre et Marie Curie (Paris 6) \\
4 Place Jussieu, \\
75005 Paris, France \\
\texttt{e-mail: \small friedland@math.jussieu.fr} \\

\medskip
\noindent
Olivier Gu\'{e}don,\\
Universit\'{e} Paris-Est Marne-La-Vall\'ee\\
Laboratoire d'Analyse et Math\'{e}matiques Appliqu\'ees, \\
5, boulevard Descartes,
Champs sur Marne, \\
77454 Marne-la-Vall\'{e}e, Cedex 2, France \\
\texttt{e-mail: \small olivier.guedon@univ-mlv.fr}
}

\newtheorem{thm}{Theorem}
\newtheorem{lemma}[thm]{Lemma}
\newtheorem{cor}[thm]{Corollary}
\newtheorem{prp}[thm]{Proposition}

\newcommand{\p}{\mathbb{P}}
\newcommand{\N}{\mathbb{N}}

\newcommand{\E}{\mathbb{E}}
\newcommand{\R}{\mathbb{R}}
\newcommand{\Id}{{\rm Id}}
\newcommand{\supp}{{\rm supp}}

\newcommand{\sparse}{{\rm sparse}}

\newcommand{\al}{\alpha}

\newcommand{\sph}{S^{n-1}}
\newcommand{\proof }{\noindent {\bf Proof. }}
\newcommand{\ST}{{\rm ST}}
\newcommand{\ess}{\rm ess}
\def \endproof {{\mbox{}\nolinebreak\hfill\rule{2mm}{2mm}\par\medbreak}}

\begin{document}

\maketitle

\footnotetext[1]{Primary Classification 
. Secondary Classification 
.}
\footnotetext[2]{Keywords. 
}

\begin{abstract}
We present a relation between sparsity and non-Euclidean isomorphic embeddings. We introduce a general restricted isomorphism property and show how it enables to construct embeddings of $\ell_p^n$, $p > 0$, into various type of Banach or quasi-Banach spaces. In particular, for $0 <r < p<2$ with $r \le 1$, we construct a family of operators that embed $\ell_p^n$ into $\ell_r^{(1+\eta)n}$, with optimal polynomial bounds in $\eta >0$.
\end{abstract}

\section{Introduction}

A quasi-Banach space $(X, \| \cdot \|)$ is said to be an $r$-normed quasi-Banach space for some $0 < r \le 1$ if: $\|x\| = 0$ iff $x=0$, $\| \lambda x \|= | \lambda| \|x\|$ for any $x \in X$, $\lambda \in \R$, and for any $x$, $y \in X$, $\| x + y \|^r \le \|x\|^r + \|y\|^r$. It is well-known \cite{Rol} that any quasi-Banach space can be equipped with an equivalent $r$-norm for a certain $r \in (0,1]$. We denote by $\sparse(m) = \{ x \in \R^n : |\supp(x)| \le m \}$ the set of vectors in $\R^n$ of cardinality of the support smaller than $m$. For $r$-normed quasi-Banach spaces $(E_1, \| \cdot \|_{E_1})$ and $(E_2, \| \cdot \|_{E_2})$ and for $ p> 0$, we define two properties of operators from $\ell_p^n$ into $E_j$, $j=1,2$, which play an important role in this paper.

We say that an operator $A : \ell_p^n \to E_1$ satisfies property ${\cal P}_1(m)$ if
$$
\forall x \in \sparse(m) \quad \al |x|_p \le \|A x\|_{E_1} \le \beta |x|_p
$$
where $ |x|_p = \left( \sum_{i=1}^n |x_i|^p \right)^{1/p}$. This property is a generalization of the Restricted Isometry Property of order $\delta$ introduced in \cite{CT}, for the Euclidean case, that is, $p=2$, $E_1 = \ell_2$, and the isometry refers to the fact that $\exists \delta \in (0,1)$ such that $\alpha = 1-\delta$ and $\beta = 1+ \delta$. We call property ${\cal P}_1(m)$ the \emph{restricted isomorphism property}. In the case $p=2$, some other versions of this property have been considered in the literature \cite{FL1, FL2}, introducing the relevance of working with general $\alpha,\beta$ and also with $E_1$ being $\ell_1^N$ instead of a Euclidean space. Here we introduce a general setting that is useful when a quasi-Banach space $E_1$ has stable type $p$ (see Section \ref{subsec:P_1} for the definition). The main difficulty is to find operators that satisfy property ${\cal P}_1(m)$ for a large $m$, and $\beta / \alpha$ being universal constant. To do so, we use random methods going back to \cite{JS,P}. For example, let $0 < r < p < 2$ with $r \le 1$, and let $E_1$ be $\ell_r^{\eta n}$ with $\eta \in (0,1]$, we exhibit families of random operators $T : \ell_p^n \to \ell_r^{\eta n}$ that satisfy property ${\cal P}_1(m)$, with overwhelming probability, for
\begin{align} \label{eq:m}
m = c_{p,r} \, \frac{\eta}{\log\left(1 + \frac{1}{\eta}\right)} \ n
\end{align}
where $c_{p,r}$ and $\beta /  \alpha$ are constants depending on $p$ and $r$. It works also in a more general setting of quasi-Banach spaces of stable type $p$.

For the second property we need the following notation. Let $x \in \R^n$ and let $\varphi_x : [n] \to [n]$ be a bijective mapping associated to a non-increasing rearrangement of $(|x_i|)$, i.e. $|x_{\varphi(1)}| \ge |x_{\varphi(2)}| \ge \cdots\ge |x_{\varphi(n)}|$. Denote by $I_k = \varphi_x(\{(k-1)m+1,\dots,km\})$ the subset of indices of the $k^{th}$ largest block of $m$ coordinates of $(|x_i|)$, for $1 \le k \le M$, where $M = \left[\frac{n}{m}\right] \le \frac{n}{m} + 1$ (note that $I_M$ may be of cardinality less than $m$). We denote by $x_{I_k}$ the restriction of $x$ to $I_k$. Clearly, $x_{I_k} \in \sparse(m)$ for $1 \le k \le M$ and
\begin{align} \label{eq:decomposition}
x = \sum_{k=1}^M x_{I_k}
\end{align}
as a disjoint sum.

We say that an operator $B : \ell_p^n \to E_2$ satisfies property ${\cal P}_2(\kappa, m)$ if
$$
\forall x \in \R^n, \quad \bigg( \sum_{k \ge 2} |x_{I_k}|_p^r \bigg)^{1/r} \le \|B x\|_{E_2} \le \left(\kappa n\right)^{1/q} |x|_p
$$
where $\frac{1}{p} + \frac{1}{q} = \frac{1}{r}$. This property with the right choice of parameters is just a simple consequence of linear algebra, it asks about finding a nice family of vectors in $E_2$. Our main simple example is that $\frac{{\rm Id_n}}{m^{1/q}} : \ell_p^n \to \ell_r^n$ satisfies property ${\cal P}_2(1/m, m)$. This is inspired by the techniques used in compressed sensing theory, see for example \cite{D, CRT}.

Now, we present our main theorem, it is a deterministic statement about Kashin-type isomorphic embedding for operators that satisfy properties ${\cal P}_1(m), {\cal P}_2(\kappa, m)$. It provides a new framework for constructing operators from $\ell_p^n$ into the quasi-Banach space $E_1 \oplus_1 E_2$, equipped with the quasi-norm $\|x\| = \|x_1\|_{E_1} + \|x_2\|_{E_2}$, where $x$ is uniquely defined by $x_1 + x_2$, $x_1 \in E_1$, $x_2 \in E_2$.

\begin{thm} \label{thm:kashin}
Let $0 < r \le p < \infty$, with $r \le 1$, and let $E_1$, $E_2$ be $r$-normed quasi-Banach spaces. Let $A : \ell_p^n \to E_1$ be an operator that satisfies property ${\cal P}_1(m)$, and let $B : \ell_p^n \to E_2$ be an operator that satisfies property ${\cal P}_2(\kappa, m)$. Denote $U = \frac{1}{\beta} \left(\frac{m}{n}\right)^{1/q} A$ and $V=\frac{1}{(\kappa n)^{1/q}} B$. Then for any $x \in \R^n$
$$
4^{-1/r} \left( \frac{\al}{\beta} \right) \left( \frac{\min(m,1/\kappa)}{n} \right)^{1/q} |x|_p \le \|Ux\|_{E_1} + \|Vx\|_{E_2} \le 3 |x|_p
$$
where $\frac{1}{p} + \frac{1}{q} = \frac{1}{r}$.
\end{thm}

As we said, we know several important situations, where we can find operators that satisfy the main properties. Let $\eta \in (0,1]$ and $Y$ be the random vector taking the values $\{\pm e_{1}, \ldots, \pm e_{\eta n}\}$, the vectors of the canonical basis in $\R^{\eta n}$, with probability $\frac{1}{2 \eta n}$. Let $(Y_{ij})$ be a sequence of independent copies of $Y$, where $1 \le i \le n$, $j \in \N$. We define the operator (see \cite{P, FG})
\begin{align} \nonumber
S : \ell_p^n & \to \ell_r^{\eta n}
\\
\label{operator:main}
x = (x_1, \ldots, x_n) & \mapsto \sum_{i=1}^n x_i \sum_{j \ge 1} \frac{1}{j^{1/p}} Y_{ij}
\end{align}

We shall prove in Section \ref{subsec:P_1} that a certain multiple of $S$ satisfies property ${\cal P}_1(m)$ with $m$ as in $(\ref{eq:m})$. An important consequence of Theorem \ref{thm:kashin} is the following:

\begin{thm} \label{thm:poly}
Let $0 < r < p < 2$ with $r \le 1$. For any $\eta \in (0,1]$ and any natural number $n$, let $W$ be a $(1+\eta) n \times n$ matrix defined by
$$
W =
\frac{1}{n^{1/q}}
\left(
\begin{array}{c}
{\rm Id_n}
\\
\tilde S
\end{array}
\right) : \ell_p^n \to \ell_r^{(1+ \eta) n}
$$
where $\tilde S = \frac{c'(p,r)}{(\log(1+1/\eta))^{1/q}} S$. Then, with probability greater than $1 - 2 \exp (-b_{p,r}\eta n)$, for any $x \in \R^n$
$$
c_{p,r} \, \left( \frac{ \eta}{\log(1+ \frac{1}{\eta})} \right)^{1/q} |x|_p \le |W x|_r \le 3 \cdot 2^{1/r} |x|_p
$$
where $\frac{1}{p} + \frac{1}{q} = \frac{1}{r}$, and $c'(p,r), b_{p,r}, c_{p,r}$ are positive constants depending on $p,r$.
\end{thm}

This answers a long standing question, whether one can give an explicit construction of random operator that embeds $\ell_p^n$ into $\ell_r^N$, where $0 < r < p < 2$, with $r \le 1$, and $N = (1+\eta)n$, $\eta \in (0,1]$, with optimal polynomial bound in $\eta$ (up to a log factor). This question was solved recently in \cite{FG}, where the isomorphism constant is $c_{p,r}^{1/\eta}$, using a full random operator (similar to the operator $S$). The improvement here comes from a reduction of the level of randomness of the operator. In a sense, it is a mixture of deterministic and random methods, which enable us to reach the best bound in the isomorphism constant. Several previous works \cite{JS, P, BB, NZ, JS2, FG} dealt with this subject. We refer to \cite{JS2, FG} for more precise references. An important remark is that the conclusion of Theorem \ref{thm:poly} holds for a lot of new operators. For example the random operators defined originally in \cite{JS} also satisfy property $P_1(m)$ with the same $m$ as in $(\ref{eq:m})$. And several other operators $B : \ell_p^n \to \ell_r^n$ satisfy property ${\cal P}_2(1/m, m)$. Hence, the strategy that we have developed allows to define several new explicit random operators that satisfy the desired conclusion.

The paper is organized as follows. In Section \ref{sec:pre}, we present the main consequence of the properties ${\cal P}_1$ and ${\cal P}_2$, that is, Theorem \ref{thm:kashin}. Of course, the delicate point is to describe some operators that satisfy the properties ${\cal P}_1(m)$ and ${\cal P}_2(\kappa, m)$ with the good parameters. This is the purpose of Corollary \ref{clm:id} in Section \ref{subsec:P_2} and Theorem \ref{thm:stable} in Section \ref{subsec:P_1}. In Section \ref{subsec:embedding}, we present the proof of Theorem \ref{thm:poly}. Finally, in the Appendix, we discuss the consequences of property ${\cal P}_1(m)$ in approximation theory and compressed sensing as it is now understood after the papers \cite{D},  \cite{CT},  \cite{CRT} and \cite{KT} (see e.g. Chapter 2 in \cite{Lamabook}). In particular, we observe that these operators are good sensing matrices when using the $\ell_r$-minimization method and that the kernel of these operators attain the optimal known bounds for the Gelfand numbers of $\Id : \ell_r^n \to \ell_p^n$. This illustrates the tightness of the method.

\section{The main theorem} \label{sec:pre}

In this section, we prove Theorem \ref{thm:kashin}. Let $x \in \R^n$ and decompose it as it is described in the introduction, see (\ref{eq:decomposition}): $x = \sum_{k=1}^M x_{I_k}$, where $(I_k)_{k=1}^M$ is the subset of indices of the $k^{th}$ largest block of $m$ coordinates of $(|x_i|)$, and $x_{I_k}$ is the restriction of $x$ to $I_k$. Each subsets $I_k$ is of cardinality $m$ except $I_M$ whose cardinality is less than $m$. Moreover, $M = \lceil \frac{n}{m} \rceil \le \frac{n}{m} +1$.

Let us start with the upper bound. By the triangle inequality,  definition of $U$ and property ${\cal P}_1(m)$, we get that
$$
\|Ux\|_{E_1}^r = \|U\sum_{k=1}^M x_{I_k}\|_{E_1}^r \le \sum_{k=1}^M \|Ux_{I_k}\|_{E_1}^r \le \left(\frac{m}{n}\right)^{r/q} \sum_{k=1}^M |x_{I_k}|_p^r
$$
Since $r \le p$ we get by H\"older's inequality
$$
\sum_{k=1}^M |x_{I_k}|_p^r \le M^{r/q} \left(\sum_{k=1}^M |x_{I_k}|_p^p\right)^{r/p}
$$
where $1/p + 1/q = 1/r$. By definition of the $\ell_p^n$-norm and of decomposition $(\ref{eq:decomposition})$ of $x$,
$$
|x|_p^p = \sum_{k=1}^M |x_{I_k}|_p^p
$$
and we get that $\|U(x)\|_{E_1} \le 2 \, |x|_p$. By definition of $V$ and property ${\cal P}_2(\kappa, m)$, we have $\|Vx\|_{E_2} \le |x|_p$. We conclude that for any $x \in \R^n$,
$$
\|Ux\|_{E_1} + \|Vx\|_{E_2} \le 3 |x|_p
$$

As for the lower bound, we partition the sphere $\sph_p$ into two sets, such that on one set we have a lower bound for $\|Ux\|_{E_1}$, and on the other set we have a lower bound for $\|Vx\|_{E_2}$. This natural type of partitioning of the sphere was also used by Kashin \cite{K} and \cite{S,AFMS}. More precisely, for $0 < \gamma < 1$ to be defined later, we partition the sphere $\sph_p$ with respect to $\gamma$ and define
$$
\Sigma_\gamma = \left\{ x \in \sph_p : \|Vx\|_{E_2} \le \gamma \right\}
$$
Clearly by this definition, if $x \notin \Sigma_{\gamma}$ then a lower bound $\gamma$ holds for this point, i.e.
$$
\|Vx\|_{E_2} > \gamma
$$
In the other case, where $x \in\Sigma_{\gamma}$, we shall obtain a lower bound, this time for the operator $U$. By the triangle inequality
$$
\|Ux\|_{E_1}^r \ge \|Ux_{I_1}\|_{E_1}^r - \|U(x-x_{I_1})\|_{E_1}^r
$$
Now, we learn each term. By decomposition $(\ref{eq:decomposition})$ of $x$, triangle inequality, and property ${\cal P}_1(m)$
\begin{align*}
& \|U(x-x_{I_1})\|_{E_1}^r = \|U \sum_{k=2}^M x_{I_k}\|_{E_1}^r \le \sum_{k=2}^M \|Ux_{I_k}\|_{E_1}^r\le \left(\frac{m}{n}\right)^{r/q} \sum_{k=2}^M |x_{I_k}|_p^r \\
& \|Ux_{I_1}\|_{E_1}^r \ge \left(\frac{\al}{\beta}\right)^r \left(\frac{m}{n}\right)^{r/q} |x_{I_1}|_p^r
\end{align*}
For $0 < r \le p < \infty$, with $r \le 1$, the $\ell_p$-norm on $\R^n$ is an $r$-norm. Hence
$$
|x_{I_1}|_p^r = |x - \sum_{k = 2}^M x_{I_k}|_p^r \ge 1 - \sum_{k= 2}^M |x_{I_k}|_p^r
$$
Therefore,
$$
\|Ux_{I_1}\|_{E_1}^r \ge \left(\frac{\al}{\beta}\right)^r \left(\frac{m}{n}\right)^{r/q} \left( 1 - \sum_{k= 2}^M |x_{I_k}|_p^r\right)
$$
Combining all the above, and since $\beta/\al \ge 1$
\begin{align*}
\|Ux\|_{E_1}^r & \ge \left(\frac{\al}{\beta}\right)^r \left(\frac{m}{n}\right)^{r/q} \left( 1 - \sum_{k= 2}^M |x_{I_k}|_p^r\right) -
\left(\frac{m}{n}\right)^{r/q} \sum_{k=2}^M |x_{I_k}|_p^r \\
& \ge \left(\frac{\al}{\beta}\right)^r \left(\frac{m}{n}\right)^{r/q} \left( 1 - \sum_{k=2}^M |x_{I_k}|_p^r \left( 1 + \left(\frac{\beta}{\al}\right)^r
\right) \right) \\
& \ge \left(\frac{\al}{\beta}\right)^r \left(\frac{m}{n}\right)^{r/q} \left( 1 - 2 \sum_{k=2}^M |x_{I_k}|_p^r \left(\frac{\beta}{\al}\right)^r \right)
\end{align*}
By property ${\cal P}_2(\kappa, m)$ and recalling that $x \in\Sigma_{\gamma}$ we have
$$
\sum_{k=2}^M |x_{I_k}|_p^r \le \|Bx\|_{E_2}^r = (\kappa n)^{r/q} \|Vx\|_{E_2}^r \le \gamma^r (\kappa n)^{r/q}
$$
It follows
$$
\|Ux\|_{E_1}^r \ge \left(\frac{\al}{\beta}\right)^r \left(\frac{m}{n}\right)^{r/q} \left( 1 - 2 \gamma^r (\kappa n)^{r/q}
\left(\frac{\beta}{\al}\right)^r \right)
$$
We conclude that if
$$
\gamma = \frac{\al}{4^{1/r} \beta} \left(\frac{1}{\kappa n}\right)^{1/q}
$$
then for any $x \in \Sigma_\gamma$
$$
\|Ux\|_{E_1} \ge \frac{\al}{2^{1/r} \beta} \, \left(\frac{m}{n}\right)^{1/q}
$$
Recalling that for any $x \notin \Sigma_{\gamma}$ we have $\|Vx\|_{E_2} > \gamma$, it implies that for any $x \in \sph_p$
$$
\frac{\al}{4^{1/r} \beta} \, \left( \frac{\min(m,1/\kappa)}{n} \right)^{1/q} \le \|Ux\|_{E_1} + \|Vx\|_{E_2}
$$
Combining with the upper bound, it concludes the proof of Theorem \ref{thm:kashin}.
\endproof

\section{Operators satisfying property ${\cal P}_2(\kappa, m)$} \label{subsec:P_2}

Property ${\cal P}_2(\kappa, m)$ can be satisfied by many operators, probably the most natural example would be the identity operator. This is just a simple consequence of the following elementary lemma about the partitioning scheme that we described above.

\begin{lemma} \label{lem:rearrang}
Let $0 < r \le p$ and let $x \in \R^n$ be decomposed as in $(\ref{eq:decomposition})$. Then for any $j \ge 1$
$$
\left( \sum_{k = j}^{M-1} |x_{I_{k+1}}|_p^r \right)^{1/r} \le \frac{1}{m^{1/q}} |x_{(I_1 \cup \ldots \cup I_{j-1})^c}|_r \le \frac{1}{m^{1/q}} |x|_r
$$
where $q$ is defined by $\frac{1}{p} + \frac{1}{q} = \frac{1}{r}$.
\end{lemma}

\proof
Let $k \ge 1$. We have
$$
|x_{I_{k+1}}|_p^r = \left( \sum_{i \in I_{k+1}} |x_i|^p \right)^{r/p} \le |I_{k+1}|^{r/p} \cdot \max_{i \in I_{k+1}} |x_i|^r
$$
and
$$
\max_{i \in I_{k+1}} |x_i|^r \le \min_{i \in I_k} |x_i|^r \le \frac{1}{|I_k|} \sum_{i \in I_k} |x_i|^r = \frac{1}{|I_k|} |x_{I_k}|_r^r
$$
We deduce
$$
\forall k \ge 1 \quad |x_{I_{k+1}}|_p^r \le \frac{|I_{k+1}|^{r/p}}{|I_k|} |x_{I_k}|_r^r \le \frac{1}{m^{r/q}} |x_{I_k}|_r^r
$$
Summing up these inequalities for all $k \ge j$ we get
$$
\sum_{k = j}^{M-1} |x_{I_{k+1}}|_p^r \le \frac{1}{m^{r/q}} \sum_{k = j}^{M-1} |x_{I_k}|_r^r \le \frac{1}{m^{r/q}} |x_{(I_1 \cup \ldots \cup I_{j-1})^c}|_r^r \le \frac{1}{m^{r/q}} |x|_r^r
$$
which concludes the proof.
\endproof

It follows that the identity operator from $\ell_p^n$ to $\ell_r^n$, correctly normalized, satisfies property $P_2(1/m,m)$.

\begin{cor} \label{clm:id}
Let $0 < r \le p < \infty$ and $q$ be such that $1/p + 1/q = 1/r$. The operator $\frac{1}{m^{1/q}} \Id_n:\ell_p^n \to \ell_r^n$ satisfies property ${\cal P}_2(\kappa, m)$, where $\kappa = \frac{1}{m}$ and $E_2 = \ell_r^n$. More precisely for any $x \in \R^n$,
$$
\sum_{k = 2}^M |x_{I_k}|_p^r \le \frac{1}{m^{r/q}}|\Id_n x|_r^r \le \left( \frac{n}{m} \right)^{r/q} |x|_p^r
$$
\end{cor}

\proof
Take $j=1$ in the Lemma \ref{lem:rearrang} and use H\"older's inequality for the upper bound.
\endproof

\bigskip

\noindent{\bf Remark.}
1.
Property ${\cal P}_2(\frac{1}{m},m)$ holds true for matrices with non-trivial kernel, e.g. take a permutation matrix
$P_\sigma$ (where $\sigma\in S_n$) and remove $k$ (up to $m$) rows.
\\
2.
Property ${\cal P}_2(\frac{K}{m},m)$ holds true for any operator $\frac{1}{m^{1/q}}V$, where $V : \ell_p^n \to E_2$ satisfies for any $x \in \R^n$
$$
|x|_r \le \|Vx\|_{E_2} \le (Kn)^{1/q} |x|_p
$$

\section{Operators satisfying property ${\cal P}_1(m)$} \label{subsec:P_1}

In order to apply Theorem \ref{thm:kashin} we should find operators that satisfy property ${\cal P}_1(m)$. For $p=2$, the set of matrices that satisfy property ${\cal P}_1(m)$ is wide. Indeed, random matrices like e.g. Bernoulli matrices, Gaussian matrices, matrices with independent log-concave rows, satisfy this property for a large $m$ with $E_1= \ell_2^d$ or $\ell_1^d$. In the case $p \ne 2$, the situation is more delicate. We shall present a possible answer when $0<p<2$, which is based on the notion of $p$-stable random variables. A natural question consists of asking what happens for $p > 2$. In that case, our method won't work, as the notion of $p$-stable random variable is not valid anymore.

\subsection{Restricted isometry property for quasi-Banach spaces}

We need several consequences of well-known results about $p$-stable random variables. We refer the reader to Chapter 5 of the book \cite{LT} and to \cite{P, BB} for the construction of the random operator that we present here. We recall that a real-valued symmetric random variable $\theta$ is called $p$-stable for $p \in (0,2]$ if its characteristic function is as follows: for some $\sigma \ge 0$, $\E \exp(it \theta) = \exp(-\sigma |t|^p)$, for any real $t$. When $\sigma = 1$, we say that $\theta$ is standard. Stable random variables are characterized by their fundamental ``stability" property: if $(\theta_i)$ is a standard $p$-stable sequence, for any finite sequence $(\al_i)$ of real numbers, $\sum_i \al_i \theta_i$ has the same distribution as $(\sum_i |\al_i|^p )^{1/p} \theta_1$.

Let $0 < r < p < 2$, with $r \le 1$, and let $X$ be an $r$-normed quasi-Banach space. We say that $X$ is of stable type $p$ if there exists a constant $\ST_p$ such that for any finite sequence $(x_i)_i \subset X$
\begin{align} \label{eq:stable type p}
\left(\E \| \sum \theta_i x_i \|^r \right)^{1/r} \le \ST_p \, (\sum \|x_i \|^p)^{1/p}
\end{align}
where $(\theta_i)_i$ is an i.i.d sequence of standard $p$-stable random variables. We denote by $\ST_p(X)$ the smallest constant $\ST_p$, such that (\ref{eq:stable type p}) holds. An important property of $p$-stable random variables is the following stability result. Let $\Theta = \sum_{j=1}^N \theta_j x_j$, where $x_j \in X$ and $\theta_j$ are i.i.d. $p$-stable random variables. For every integer $k \ge 1$, if $\Theta_1, \ldots, \Theta_k$ are independent copies of $\Theta$ then for every $(\al_i)_{i=1}^k \in \R^k$, $\sum_{i=1}^k \al_i \Theta_i$ has the same distribution as $(\sum_{i=1}^k |\al_i|^p)^{1/p} \, \Theta$. In particular,
\begin{align} \label{eq:stability}
\left( \E \| \sum_{i=1}^k \al_i \Theta_i \|^r \right)^{1/r} = (\sum_{i=1}^k |\al_i|^p )^{1/p} ~ \left( \E\| \Theta \|^r \right)^{1/r}
\end{align}

Assume that $X$ is of stable type $p$ with $0< r < p< 2$. Therefore, we can find a finite sequence $x_1, \ldots, x_N \in X$ such that $\sum_{i=1}^N \|x_i\|^p  = 1$, and 
\begin{align} \label{eq:typedef}
\left( \E\|\sum_{i=1}^N \theta_i x_i \|^r \right)^{1/r} \ge \frac{1}{2}\ST_p(X) 
\end{align}
Let $y_i = x_i/\|x_i\|$. Let $Y$ be a symmetric $X$-valued random vector with distribution equal to $\sum_{i=1}^N \|x_i\|^p (\delta_{y_i} + \delta_{-y_i})/2$, and let $Y_1, Y_2, \ldots$ be i.i.d copies of $Y$.
Let $(\lambda_i)$ be independent random variables with common exponential distribution $\p \{ \lambda_i > t\} = \exp(-t)$, $t\ge 0$. Set $\Gamma_j = \sum_{i=1}^j \lambda_i$, for $j \ge 1$
then it is known (cf. \cite{LWZ}) that there exists a positive constant $s'_p$ depending on $p$, such that in distribution
$$
\tilde{\Theta} = \sum_{j \ge 1} \Gamma_j^{-1/p} Y_j {\overset{d}{~=~}} s'_p \sum_{i=1}^N \theta_i x_i
$$
It follows that
\begin{align} \label{eq:position}
\left(\E \| \tilde{\Theta} \|^r \right)^{1/r} \ge s_p \cdot \ST_p(X)
\end{align}

Following Pisier \cite{P}, we define the operator
\begin{align} \nonumber
T : \ell_p^n & \to X
\\
\label{operator}
\al = (\al_1, \ldots, \al_n) & \mapsto \frac{1}{\left(\E \| \tilde{\Theta} \|^r \right)^{1/r}} \sum_{i=1}^n \al_i \sum_{j \ge 1} j^{-1/p} Y_{ij}
\end{align}

\begin{thm} \label{thm:stable}
Let $0 < \frac{2p}{p+2} < r < p < 2$, with $r \le 1$, and let $X$ be an $r$-normed quasi-Banach space, with stable type $p$ constant $\ST_p(X)$. Then with probability greater than $1 - 2 \exp \left( -b_{p,r} \left( \ST_p(X) \right)^q \right)$ the operator $T$ satisfies property ${\cal P}_1(m)$ for any $m$, such that
$$
m \le {\left( c_{p,r} \ST_p(X)\right)^q} \, / \ {\log\left( 1 + {n}/{\left( c_{p,r} \ST_p(X)\right)^q} \right)}
$$
where $\frac{1}{p} + \frac{1}{q} = \frac{1}{r}$ and $b_{p,r}, c_{p,r}$ are positive constants depending on $p$ and $r$. More precisely, for such $m$
$$
\forall \al \in \sparse(m) \quad \frac{1}{2}\ |\al|_p^r \le \|T\al\|^r \le \frac{3}{2} |\al|_p^r
$$
\end{thm}

\bigskip

\noindent{\bf Remark.}
\\
1. For any $\delta \in (0,1)$, we can introduce a dependence in $\delta$ in the choice of $m$, such that property ${\cal P}_1(m)$ holds with $\al = 1-\delta$ and $\beta = 1+ \delta$. So, this is an extension of the RIP to $r$-normed quasi-Banach spaces.
\\
2. Notice that $(\ref{eq:typedef})$ is the main property that should be satisfied by the family of vectors $\{ x_1, \ldots, x_N \} \subset X$. The quantity $\ST_p(X)$ in the theorem can be replaced by the quantity that will appear in this inequality, for the prescribed family $\{ x_1, \ldots, x_N \}$. It is not necessary to estimate $\ST_p(X)$ but  the definition of stable type $p$ corresponds to make it as large as possible.

\bigskip

For the proof of Theorem \ref{thm:stable}, we define the following auxiliary operator
\begin{align*}
\tilde{T} : \ell_p^n & \to X \\
\al = (\al_1, \ldots, \al_n) & \mapsto \frac{1}{\left(\E \| \tilde{\Theta} \|^r \right)^{1/r}} \sum_{i=1}^n \al_i \sum_{j \ge 1} \Gamma_j^{-1/p} Y_{ij}
\end{align*}

We need the following two lemmas, which are analogous to the main lemmas in \cite{P, BB}. The first one is a consequence of well-known results about $p$-stable random variables (see also \cite[Lemma 0.3]{BB}).

\begin{lemma} \label{lem:comparison}
Let $0 < \frac{2p}{p+2} < r < p < 2$. For the operators $T$ and $\tilde T$ defined above, there exists a positive constant $d_{p,r}$ depending on $p$ and $r$, such that for any $\al \in \R^n$,
$$
\left| \E \|T\al \|^r - \E \| \tilde{T} \al \|^r \right| \le
\frac{d_{p,r} }{\E \| \tilde{\Theta} \|^r} \sum_{i=1}^n |\al_i|^r
\le
\frac{d_{p,r} ~ |\supp(\al)|^{r/q}}{\E \| \tilde{\Theta} \|^r}
$$
\end{lemma}

\proof
It is known (see \cite{FG} for details) that for any $0 < \frac{2p}{p+2} < r < p < 2$ there exists a positive constant $d_{p,r}$ depending on $p$ and $r$, such that
$$
\sum_{j \ge 1} \E \left|j^{-1/p} - \Gamma_j^{-1/p} \right|^r \le d_{p,r}
$$
Therefore,
\begin{align*}
\left| \E \|T\al \|^r - \E \| \tilde{T} \al \|^r \right| & \le \frac{1}{\E \| \tilde{\Theta} \|^r} \E \left\| \sum_{i \in \supp(\al)} \al_i \sum_{j \ge 1} \left( j^{-1/p}Y_{ij} - \Gamma_j^{-1/p} Y_{ij}\right) \right\|^r \\
& \le \frac{1}{\E \| \tilde{\Theta} \|^r} \sum_{i \in \supp(\al)} |\al_i|^r ~ \E \left\| \sum_{j \ge 1} \left( j^{-1/p}Y_{ij} - \Gamma_j^{-1/p}Y_{ij} \right) \right\|^r \\
& \le \frac{1}{\E \| \tilde{\Theta} \|^r} \sum_{i \in \supp(\al)} |\al_i|^r \sum_{j \ge 1} \E \left| j^{-1/p} - \Gamma_j^{-1/p} \right|^r \|Y_{ij}\|^r \\
& \le \frac{1}{\E \| \tilde{\Theta} \|^r} \sum_{i \in \supp(\al)} |\al_i|^r d_{p,r} \le \frac{d_{p,r} ~ |\supp(\al)|^{r/q}}{\E \| \tilde{\Theta} \|^r}
\end{align*}
\endproof

The next lemma follows from results about scalar martingale difference (cf. \cite{JS,P, BB}).

\begin{lemma} \label{lem:concentration}
Let $X$ be an $r$-normed quasi-Banach space and $(Z_j)_j$ be a sequence of independent $X$ valued random vectors, which are uniformly bounded. Let $\lambda_k = \ess \sup \|Z_k (\cdot)\|$. If $Z = \sum_{k \ge 1} Z_k$ converges a.s. then for any $t>0$ we have
$$
\p \left\{ \left| \|Z\|^r - \E\|Z\|^r \right| \ge t \right\} \le 2 \exp \left( -c'_{p,r} \left( \frac{t}{\|(\lambda_k^r)_k\|_{p/r,\infty}} \right)^{q/r} \right)
$$
where $c'_{p,r}$ is a positive constant depending on $p$ and $r$.
\end{lemma}

Denote $Z_k = \al_i j^{-1/p} Y_{ij}$ and observe that $\|Z_k\| \le |\al_i| j^{-1/p}$. It is easy (cf. \cite{P,FG}) to deduce that $\|(\lambda_k^r)_k\|_{p/r,\infty} \le |\al|_p = 1$. Therefore, by applying Lemma \ref{lem:concentration}, for any $t>0$, we obtain
$$
\p \left\{ \left| \|T\al\|^r - \E\|T\al\|^r \right| \ge \frac{t}{\E \| \tilde{\Theta} \|^r} \right\} \le 2 \exp(-c'_{p,r} t^{q/r})
$$
Taking $t = s^r ~ \E\|\tilde{\Theta}\|^r$, we conclude that for any $\al \in S_p^{n-1}$,
\begin{align} \label{eq:concentration}
\p \left\{ \left| \|T\al\|^r - \E\|T\al\| \right|^r \ge s^r \right\} \le 2 \exp \left( -c'_{p,r} \, s^q \left( \E\|\tilde{\Theta}\|^r \right)^{q/r} \right)
\end{align}

\noindent {\bf Proof of Theorem \ref{thm:stable}.}
Let $\al \in S_p^{n-1}$. By (\ref{eq:stability}), and the discussion above we have
$$
\E \| \tilde{T} \al \|^r = \frac{1}{\E \| \tilde{\Theta} \|^r} \E \| \sum_{i=1}^n \al_i \sum_{j \ge 1} \Gamma_j^{-1/p} Y_{ij} \|^r = \frac{1}{\E \| \tilde{\Theta} \|^r} \E \| \sum_{i=1}^n \al_i \tilde{\Theta_i} \|^r = \frac{\E \| \tilde{\Theta} \|^r} {\E\|\tilde{\Theta}\|^r} = 1
$$
Therefore, by Lemma \ref{lem:comparison}
$$
\left| \E\|T\al\|^r - 1 \right| = \left| \E \|T\al \|^r - \E \| \tilde{T} \al \|^r \right| \le \frac{d_{p,r} |\supp(\al)|^{r/q}}{\E\|\tilde{\Theta}\|^r}
$$
It follows that if $|\supp(\al)| \le \left( \frac{\E\|\tilde{\Theta}\|^r}{4 d_{p,r}} \right)^{q/r}$ then
$$
\left| \E\|T\al\|^r - 1 \right| \le \frac{1}{4}
$$
Moreover, by (\ref{eq:concentration}) for $s=(1/8)^{1/r}$
$$
\p \left\{ \left| \|T\al\|^r - \E\|T\al\|^r \right| \ge 1/8 \right\} \le 2 \exp \left( -b'_{p,r} \left( \E\|\tilde{\Theta}\|^r \right)^{q/r} \right)
$$
We deduce that for every $\al \in S_p^{n-1}$, such that $|\supp(\al)| \le \left( \frac{\E\|\tilde{\Theta}\|^r}{4 d_{p,r}} \right)^{q/r}$
$$
\p \left\{ 5/8 \le \|T\al\|^r \le 11/8 \right\} \ge 1 - 2 \exp \left( -b'_{p,r} \left( \E\|\tilde{\Theta}\|^r \right)^{q/r} \right)
$$
We need to approximate the set of sparse vectors of size $m$ of $S_p^{n-1}$ by a net. By a $\delta$-net of a subset $U$ of an $r$-normed space $X$, we mean a subset $\cal N$ of $U$, such that for all $x \in U$,
$$
\inf_{y \in {\cal N}} \|x-y\|^r \le \delta
$$
It is well-known by a volumetric argument (see \cite[Lemma 2]{JS}) that if $X$ is an $r$-normed space of dimension $m$ then the unit sphere of $X$ contains a $\delta$-net of cardinality at most $(1+2/\delta)^{m/r}$. Now, since $$
\sparse(m) \cap S_p^{n-1} = \bigcup_{| I |= m} \R^I \cap S_p^{n-1}
$$
we can find a $1/12$-net $\cal N$ of $\sparse(m) \cap S_p^{n-1}$ of the form $\cup_{| I |= m} {\cal N}_I$, where ${\cal N}_I$ is a subset of $\R^I \cap S_p^{n-1}$ of cardinality at most $25^{m/r}$. The cardinality of $\cal N$ is at most $\binom{n}{m} 25^{m/r} \le \exp(m \log(d_r m/n))$, where $d_r$ is a constant depending on $r$. From a classical union bound argument, we deduce that with probability greater than
$$
1 - 2 \exp \left( -b'_{p,r} \left( \E\|\tilde{\Theta}\|^r \right)^{q/r} + m \log\left(\frac{d_r m}{n}\right) \right)
$$
we have
$$
\forall y \in {\cal N} , \, \, 5/8 \le \|T y\|^r \le 11/8
$$
Since for any $\al \in \sparse(m) \cap S_p^{n-1}$ there exists $y \in {\cal N}$, such that
$\al-y \in \sparse(m)$ and $|\al-y|_p^r \le 1/12$. And, for $r \le p$, the $\ell_p$-norm is an $r$-norm. We get by the triangle inequality of the $r$-norm $\|\cdot\|$
\begin{align*}
\forall \al \in \sparse(m) \cap S_p^{n-1}, \|T \al\|^r \ge \| T y \|^r - \|T(\al-y)\|^r \\
{\rm and} \ \sup_{\al \in \sparse(m) \cap S_p^{n-1}} \|T \al \|^r \le
\sup_{y \in {\cal N}} \|T y\|^r + \frac{1}{12} \sup_{\al \in \sparse(m) \cap S_p^{n-1}} \|T \al \|^r
\end{align*}
It is easy to conclude that with probability greater than
$$
1 - 2 \exp \left( -b'_{p,r} \left( \E\|\tilde{\Theta}\|^r \right)^{q/r} + m \log\left(\frac{d_r m}{n}\right) \right)
$$
we have
$$
\forall \al \in \sparse(m) \cap S_p^{n-1} , \, \, 1/2 \le \|T\al\|^r \le 3/2
$$
By $ (\ref{eq:position})$ we know that
$$
\left( \E\|\tilde{\Theta}\|^r \right)^{1/r} \ge s_p \ST_p(X)
$$
and we conclude that for constants $b_{p,r}$ and $c_{p,r}$ depending on $p$ and $r$, if
$$
m \le \left( c_{p,r} \ST_p(X) \right)^q \, / \ {\log\left( 1 + {n}/{\left( c_{p,r} \ST_p(X)\right)^q} \right)}
$$
then
$$
\p \left\{ \forall \al \in \sparse(m) \cap \sph_p , \, 1/2 \le \|T\al\|^r \le 3/2 \right\} \ge 1 - 2 \exp \left( - b_{p,r} \left( \ST_p(X) \right)^q \right)
$$
This ends the proof.
\endproof

\subsection{Restricted isomorphism property for $\ell_r$}

Let $0 < r < p < 2$ with $r \le 1$, let $\eta \in (0,1]$ and $X=\ell_r^{\eta n}$. It's well-known that $\ST_p(\ell_r^{\eta n}) = c'_{p,r} (\eta n)^{1/q}$. It's easy to see from definition $(\ref{eq:stable type p})$ that we may take the canonical basis of $\R^{\eta n}$ as the $x_i$'s in $(\ref{eq:typedef})$. Hence, let $Y$ be the random vector taking the values $\{\pm e_{1}, \ldots, \pm e_{\eta n}\}$, the vectors of the canonical basis in $\R^{\eta n}$, with probability $\frac{1}{2 \eta n}$. We define the operator (see also \cite{FG})
\begin{align*}
S : \ell_p^n & \to \ell_r^{\eta n}
\\
x = (x_1, \ldots, x_n) & \mapsto \sum_{i=1}^n x_i \sum_{j \ge 1} \frac{1}{j^{1/p}} Y_{ij}
\end{align*}

We deduce from Theorem \ref{thm:stable} the following important corollary.

\begin{cor} \label{cor:rip}
Let $0 < r < p < 2$, with $r \le 1$. Let $\eta \in (0,1]$ and $\delta = c_{p,r} \eta / \log(1+{1}/{\eta})$. Then with probability greater than $1 - 2 \exp(-b_{p,r} \eta n)$, the operator $S / (\eta n)^{1/q}$ satisfies property ${\cal P}_1(\delta n)$. More precisely,
$$
\forall x \in \sparse(\delta n) \quad c(p,r) |x|_p \le \frac{1}{(\eta n)^{1/q}} |Sx|_r \le C(p,r) |x|_p
$$
where $c_{p,r}, c(p,r), C(p,r)$ and $b_{p,r}$ are positive constants depending on $p$ and $r$.
\end{cor}

\proof
It is important to note that the definition of $S$ does not depend on the choice of $r$. Let $r \in (0,1]$.
If $0 < \frac{2p}{p+2} < r < p < 2$, this is a direct application of Theorem \ref{thm:stable} and the fact that
$\ST_p(\ell_r^{\eta n}) = c'_{p,r} (\eta n)^{1/q}$.
For the other values of $r$, we use a classical extrapolation trick. Let $r_1 \le 1$ and $r_2 \le 1$ be such that $0 < \frac{2p}{p+2} < r_1 < r_2 < p < 2$ then we can use the first case and deduce that with probability greater than $1 - 4 \exp(-b_{p} \eta n)$
$$
\forall \al \in \sparse(\delta n)
\, \
\left\{
\begin{array}{l}
\displaystyle
c_1(p) |\al|_p \le \frac{1}{(\eta n)^{1/q_1}} |S \al|_{r_1} \le C_1(p) |\al|_p
\\
\\
\displaystyle
c_2(p) |\al|_p \le \frac{1}{(\eta n)^{1/q_2}} |S \al|_{r_2} \le C_2(p) |\al|_p
\end{array}
\right.
$$
where $\frac{1}{p} + \frac{1}{q_1} = \frac{1}{r_1}$, $\frac{1}{p} + \frac{1}{q_2} = \frac{1}{r_2}$ and $b_p = \min(b_{p,r_1}, b_{p,r_2})$. For any $r<r_1$, we have for any $z \in \R^{\eta n}$, $|z|_{r_1} \le |z|_r^\theta |z|_{r_2}^{1-\theta}$, where $1/r_1 = \theta/r + (1-\theta)/r_2$, and $|z|_r \le (\eta n)^{1/r-1/r_1} |z|_{r_1}$. It is easy to deduce from the previous inequalities that
$$
\forall \al \in \sparse(\delta n), \ c(p)^{1/ r} \, |\al|_p \le \frac{1}{(\eta n)^{1/q}} |S \al|_{r} \le C_1(p) \, |\al |_p
$$
for a new positive number $c(p)$ depending on $p$.
\endproof

\medskip

\noindent {\bf Remark.} We should add that the random operator defined in \cite{JS} also satisfy the same property ${\cal P}_1(m)$, since the main properties of this random operator are completely analogous to Lemmas \ref{lem:comparison} and \ref{lem:concentration}.

\section{Random embedding of $\ell_p^n$ into $\ell_r^N$} \label{subsec:embedding}

In this section, we describe the main consequences of properties ${\cal P}_1$ and ${\cal P}_2$ in the geometry of Banach spaces. We prove Theorem \ref{thm:poly}, about the existence of very tight embeddings from $\ell_p^n$ into $\ell_r^{(1+\eta) n}$, where $\eta$ is an arbitrary small number, and $0<r < p <2$ with $r \le 1$. The isomorphism constant of these operators is of the order of $(\eta / \log(1+ 1/ \eta))^{-1/q}$, where $1/p + 1/q = 1/r$. Up to the logarithmic term it is best possible since the Banach Mazur distance between $\ell_p^n$ and $\ell_r^{n+1}$ is of the order $n^{1/q}$. It improves the main result of \cite{FG} and we refer to this paper and to \cite{JS2} for the history of the problem.

Let $S : \ell_p^n \to \ell_r^{\eta n}$ be the operator defined in $(\ref{operator:main})$ and $W : \ell_p^n \to \ell_r^{(1+\eta) n}$ be defined by
$$
W =
\frac{1}{n^{1/q}}
\left(
\begin{array}{c}
{\rm Id_n}
\\
\tilde S
\end{array}
\right) : \ell_p^n \to \ell_r^{(1+ \eta) n}
$$
where $\tilde S = \frac{c'(p,r)}{(\log(1+1/\eta))^{1/q}} S$ and $c'(p,r)$ is a constant depending on $p$ and $r$. By Corollary \ref{cor:rip}, with probability greater than $1 - 2 \exp(-b_{p,r} \eta n)$, the operator $S/(\eta n)^{1/q}$ satisfies property ${\cal P}_1(\delta n)$, that is
$$
\forall x \in \sparse(\delta n) \quad c(p,r) |x|_p \le \frac{1}{(\eta n)^{1/q}} |S x|_r \le C(p,r) |x|_p
$$
where $\delta = c_{p,r} \eta / \log(1+ {1}/{\eta})$ and $b_{p,r}, c_{p,r}, c(p,r), C(p,r)$ are numbers depending on $p$ and $r$. Moreover, by Corollary \ref{clm:id}, the operator $\frac{1}{(\delta n )^{1/q}}{\rm Id_n} : \ell_p^n \to \ell_r^n$ satisfies property ${\cal P}_2(1/\delta n, \delta n)$. Therefore, by Theorem \ref{thm:kashin}, we conclude that, for any $x \in \R^n$
$$
4^{-1/r} \frac{c(p,r)}{C(p,r)} \, \left( \frac{ \eta}{\log(1+ \frac{1}{\eta})} \right)^{1/q} |x|_p \le \frac{1}{n^{1/q}} \left( |x|_r + |\tilde Sx|_r \right) \le 3 |x|_p
$$
for $\frac{c'(p,r)}{(\log(1+1/\eta))^{1/q}} S = \tilde S$. Since
$$
|x|_r + |\tilde Sx|_r \le \left( |x|_r^r + |\tilde Sx|_r^r \right)^{1/r} = n^{1/q} |Wx|_r \le 2^{1/r} \left( |x|_r + |\tilde Sx|_r \right)
$$
the proof of Theorem \ref{thm:poly} is complete.
\endproof

\section*{Appendix} \label{sec:gelfand}

We would like to conclude this paper by presenting some relations between property ${\cal P}_1$ and between approximation theory and compressed sensing. Indeed, in his seminal paper \cite{D}, Donoho made several connections between compressed sensing and Gelfand numbers. We pursue that direction and present direct consequences of property ${\cal P}_1$ in this setting, like in \cite{D, CT, FL2, FPRU}. The results are known. The main point is to emphasize about new subspaces that "attain" Gelfand widhts, and new operators that satisfy the approximate reconstruction property via the $\ell_r$-minimization method, for $0<r \le 1$.

Recall that the Gelfand numbers of an operator $u:X \to Y$ is defined for every $k \in \N$ by
$$
c_k(u: X \to Y) = \inf\{ \sup_{x \in S, \|x\|_X \le 1} \|u(x)\|_Y\}
$$
where the infimum runs over all subspaces $S$ of codimension striclty less than $k$. For any $0<r\le \infty$, we define the weak-$\ell_r^n$ space to be $\R^n$ equipped with the quasi-norm $|\cdot|_{r, \infty}$
$$
\forall x \in \R^n, \quad |x|_{r, \infty} = \max_{i=1, \ldots, n} i^{1/r} x_i^*
$$
where $x_1^* \ge x_2^* \ge \dots \ge x_n^*$ is the non-increasing rearrangement of $(|x_i|)_{i=1}^n$.

\begin{prp} \label{prop:nullspace}
Let $E_1$ be an $r$-normed quasi-Banach space, $A : \ell_p^n \to E_1$ with $0 < r \le p$ and $m \le n$, such that property ${\cal P}_1(m)$ holds true. Then
\begin{align} \label{eq:nullspace1}
\forall \, h \in \ker A, \quad |h|_p
\le
 \left( 1 + (\beta/\al)^p \right)^{1/p} \left( \sum_{k = 2}^M |h_{I_k}|_p^r \right)^{1/r}
\end{align}
where $h = \sum_{k =1}^M h_{I_k}$ is the decomposition defined in $(\ref{eq:decomposition})$ and $M = \lceil \frac{n}{m} \rceil$.
Let $s \le m$ then
\begin{align} \label{eq:nullspace2}
\forall \, h \in \ker A, \forall \, I \subset \{1, \ldots, n\}, | I | \le s, \ |h_I|_r \le \left( \frac{s}{m} \right)^{1/q} \left( 1 + (\beta/\al)^p \right)^{1/p} |h|_r
\end{align}
\end{prp}

\noindent
{\bf Remark.} The conclusion of the proposition does not depend on the choice of $E_1$. It is chosen such that the property ${\cal P}_1(m)$ holds true for  a value of $m$ as large as possible.

\proof
Let $h \in \ker A$, $h \not= 0$ and write
$$
|h |_p^p = |h - h_{I_1}|_p^p + |h_{I_1}|_p^p
$$
Since $r \le p$, the first term satisfies
$$
|h - h_{I_1}|_p = \left|\sum_{k = 2}^M h_{I_k} \right|_p = \left( \sum_{k = 2}^M |h_{I_k}|_p^p \right)^{1/p} \le\left( \sum_{k = 2}^M |h_{I_k}|_p^r \right)^{1/r}
$$
For the second term, we use property ${\cal P}_1(m)$. Since $h \in \ker A$ then $A (h_{I_1}) = - \sum_{k=2}^{M} A(h_{I_k})$ and since $h_{I_1} \in \sparse(m)$ we get by property ${\cal P}_1(m)$
$$
|h_{I_1}|_p \le \frac{1}{\al} \| A(h_{I_1}) \|_{E_1} = \frac{1}{\al} \left\| \sum_{k=2}^{M} A(h_{I_k}) \right\|_{E_1}
\le
\frac{1}{\al} \left( \sum_{k=2}^{M} \| A(h_{I_k}) \|_{E_1}^r \right)^{1/r}
$$
where the last inequality comes from the triangle inequality for the quasi-norm in $E_1$. Since for any $k \ge 2$, $h_{I_k} \in \sparse(m)$, we have by property ${\cal P}_1(m)$, for any $k \ge 2$
$$
\| A(h_{I_k}) \|_{E_1}^r \le \beta^r |h_{I_k}|_p^r
$$
Combining both terms, It is easy to deduce $(\ref{eq:nullspace1})$. If $h \in \ker A$ and $I \subset \{1, \ldots, n\}$ with $| I | \le s$, we have by H\"older and $(\ref{eq:nullspace1})$
$$
| h_I |_r \le s^{1/q} | h |_p \le s^{1/q} \left( 1 + (\beta/\al)^p \right)^{1/p} \left( \sum_{k = 2}^M |h_{I_k}|_p^r \right)^{1/r}
$$
From Lemma \ref{lem:rearrang}, we conclude that
$$
|h_I|_r \le \left( \frac{s}{m} \right)^{1/q} \left( 1 + (\beta/\al)^p \right)^{1/p} |h|_r
$$
\endproof

The argument that we have presented goes back to  \cite{D} while working in a Euclidean setting. We refer also to \cite{Lamabook} for more details. Inequality $(\ref{eq:nullspace1})$ has some obvious consequences in terms of Gelfand numbers of the identity operator between some sequence spaces that are known from \cite{FPRU}. Inequality $(\ref{eq:nullspace2})$ is a strong form of the so called {\em null space property} and has consequences in compressed sensing.

\begin{cor}
Let $0 < r < p < 2$ with $r < 1$. Then for any $c_{p,r} \log n \le k \le n-1$,
\begin{align} \label{eq:c_k}
c_k(\Id : \ell_r^n \to \ell_p^n) \le c_k(\Id : \ell_{r, \infty}^n \to \ell_p^n) \le C_{p,r} \left( \frac{\log(1 + \frac{n}{k}) }{k} \right)^{1/q}
\end{align}
where $1/p + 1/q = 1/r$, and $c_{p,r}, C_{p,r}$ are positive constants depending on $p,r$.
\end{cor}

The upper bound in $(\ref{eq:c_k})$ is known to be optimal \cite{FPRU} (up to constants depending on $p$ and $r$), it is proved by interpolation in \cite{GueLit} (see also \cite{V}).  Here, we give an alternative proof based on our method, i.e. we find new subspaces for which this bound is attained.

\medskip

\proof
For any $x \in \R^n$, $|x|_{r, \infty} \le |x|_r$. Hence, the left inequality is obvious. Let $c_{p,r} \log n \le k \le n-1$ and $E_1 = \ell_{r_1}^{k}$, where $r < r_1 < p$ and $r_1 \le 1$. Let $\eta$ such that $k = \eta n$, we know from Corollary \ref{cor:rip} that $S/k^{1/r_1-1/p}$ satisfies property ${\cal P}_1(m)$, where $m = c_{p,r_1} k / \log(1+n/k)$, $\al$ and $\beta$ being constants depending on $p$ and $r_1$. Following the proof of Lemma \ref{lem:rearrang}, we know that for any $h \in \R^n$,
$$
| h_{I_{k+1}}|_p \le m^{1/p} \, h_{km}^* \le \frac{k^{-1/r}}{m^{1/q}} \, |h|_{r, \infty}
$$
by definition of the weak-$\ell_r^n$ norm. Since $r_1 / r > 1$, $\sum_{k \ge 1} k^{-r_1/r}$ is finite and
$$
\left( \sum_{k = 2}^M |h_{I_k}|_p^{r_1} \right)^{1/r_1} \le \frac{c_{p,r}}{m^{1/q}} \, |h|_{r, \infty}
$$
By definition of $S$, ${\rm codim} \ker S < k$ and we conclude by Proposition \ref{prop:nullspace} and $(\ref{eq:nullspace1})$ that
$$
\hbox{ for any } h \in \ker S, \quad |h|_p \le C_{p,r} \left( \frac{\log(1 + \frac{n}{k}) }{k} \right)^{1/q} |h|_{r, \infty}
$$
which ends the proof.
\endproof

\bigskip

\noindent{\bf Remark.} Obviously, for $r=1$, we get with the same proof an additional $\log(1+n/k)$ factor as in \cite{FPRU}.

\bigskip

There is a connection with the analysis of sparse recovery via $\ell_r$-minimization method for $0 < r \le 1$. Let $S : \R^n \to \R^k$ and for any $y \in \R^n$,
\begin{align} \label{eq:reconstruction}
\Delta_r(y) = \hbox{argmin } |z|_r, \hbox{ subject to } Sz = Sy
\end{align}

For $r=1$, this is the basis pursuit algorithm \cite{CRT}, and it has been generalized to the non-convex minimization problem for $r < 1$ in \cite{Chartrand, FL1, FPRU}. It is known \cite{FPRU} that if a matrix satisfies the Restricted Isometry Property then the $\ell_r$-minimization method gives good approximate reconstruction of signals (which is exact in the sparse case). Since inequality $(\ref{eq:nullspace2})$ is the analogue of the strong form of the so called {\em null space property} in compressed sensing, we conclude that an operator satisfying property ${\cal P}_1(m)$ for some quasi-Banach space $E_1$ is a good sensing matrix. We illustrate it in the following corollary.

\begin{cor}
Let $S : \ell_p^n \to \ell_r^k$ be the random operator defined in $(\ref{operator:main})$.
If $s>0$ satisfies
$$
s \le c(p,r) \frac{k}{\log(1+\frac{n}{k})}
$$
then with probability greater than $1 - \exp(-b_{p,r} k)$,
$$
\left| y - \Delta_r (y)\right|_r \le 4^{1/r} \inf_{|I| \le s} |y - y_I|_r
$$
And if $y \in \sparse(s)$, the reconstruction is exact: $y = \Delta_r(y)$.
\end{cor}

\proof
By definition of $\Delta_r$, $h = y - \Delta_r(y) \in \ker S$ and $|\Delta_r(y)|_r \le |\Delta_r(y) + h|_r$.
We get
$$
|y|_r^r \ge |y + h |_r^r = |y_I + h_I + y_{I^c} + h_{I^c}|_r^r
\ge
| y_{I}|_r^r -| h_{I}|_r^r +|h_{I^c}|_r^r-| y_{I^c}|_r^r
$$
so that
\begin{align} \label{eq:1}
2 |y_{I^c}|_r^r \ge |h_{I^c}|_r^r - | h_{I}|_r^r
\end{align}
By Corollary \ref{cor:rip}, we know that with probability greater than $1 - \exp(-b_{p,r} k)$, the operator
$S/k^{1/q}$ satisfies property ${\cal P}_1(m)$, where $m = c_{p,r} k / \log(1+n/k)$, $\al$ and $\beta$ being constants depending on $p$ and $r$.
We can apply Proposition \ref{prop:nullspace} and we deduce from $(\ref{eq:nullspace2})$
that
$$
|h_I|_r \le
\left( \frac{s}{m} \right)^{1/q}
\left( 1 + (\beta/\al)^p \right)^{1/p}
|h|_r
$$
We choose the constant $c(p,r)$ in the definition of $s$, such that
$$
\left( \frac{s}{m} \right)^{1/q}
\left( 1 + (\beta/\al)^p \right)^{1/p}
\le \frac{1}{4^{1/r}}
$$
hence, we get that $|h_I|_r^r \le |h|_r^r / 4$, which gives that $|h_I|_r^r \le |h_{I^c}|_r^r / 3$. We conclude from $(\ref{eq:1})$ that
$
|y_{I^c}|_r^r \ge |h_{I^c}|_r^r /3
$ and that
$$
|h|_r^r = |h_I|_r^r + |h_{I^c}|_r^r \le 4 |y_{I^c}|_r^r
$$
which is the announced result.
\endproof

\address

\end{document}